\documentclass[12pt]{amsart}

\usepackage[arrow,dvips,matrix,arrow,ps,color,line,curve,frame]{xy}

\usepackage{amssymb}

\advance\textheight1cm

\let\oldlabel=\label
\def\prellabel{\marginparsep=1em
    \def\label##1{\oldlabel{##1}\ifmmode\else\ifinner\else
         \marginpar{{\footnotesize\ \\ \tt
                    ##1}}\fi\fi}}

\advance\headheight1.15pt

\let\:\colon
\let\epsilon\varepsilon

\let\phi=\varphi
\let\theta=\vartheta

\let\Bbb=\mathbb

\def\opn#1#2{\def#1{\operatorname{#2}}}
\opn\gp{gp} \opn\Max{Max} \opn\Ker{Ker} \opn\Coker{Coker}
\opn\Ext{Ext} \opn\conv{conv} \opn\chara{char} \opn\n{n} \opn\h{h}
\opn\GL{GuL} \opn\SL{SL} \opn\sn{sn} \opn\inte{int} \opn\End{End}
\opn\rank{rank} \opn\Aff{Aff} \opn\Spec{Spec} \opn\Proj{Proj}
\opn\QF{QF} \opn\I{Im} \opn\Hom{Hom} \opn\Aut{Aut} \opn\W{Witt}
\opn\W{W} \opn\inte{int} \opn\pyr{pyr} \opn\l{l} \opn\r{r}
\opn\const{const}

\def\QQ{{\Bbb Q}}

\opn\End{End}
\opn\U{U}%
\opn\ch{ch}%

\def\Q{{\Box\kern1pt}}%

\def\kk{{\bf k}}

\def\F{{\textsf F}}

\newtheorem{lemma}{Lemma}
\newtheorem{corollary}[lemma]{Corollary}
\newtheorem{theorem}[lemma]{Theorem}
\newtheorem{proposition}[lemma]{Proposition}

\theoremstyle{definition}

\newtheorem{remark}[lemma]{Remark}

\begin{document}

\title[On Bass' question for finitely generated algberas]{On Bass' question for finitely generated\\ algebras
over large fields}

\author{Joseph Gubeladze}

\thanks{Supported by NSF grant DMS-0600929}


\address{Department of Mathematics, San Francisco
State University, San Francisco, CA 94132, USA}

\email{soso@math.sfsu.edu}

\begin{abstract}
Recently Corti\~nas-Haesemayer-Walker-Weibel  gave affirmative
answer to Bass' 1972 question on $NK$-groups for algebras of
essentially finite type over large fields of characteristic 0.
Here we give an alternative short proof of this result for
algebras of finite type over such fields. Our approach is based on
classical techniques in higher $K$-theory of rings and a direct
$K_i$-analogue of an old observation of Murthy-Pedrini, dating
back from the same 1972.
\end{abstract}

\maketitle

Recently, for algebras $R$ over the rationals
Corti\~nas-Haesemeyer-Walker-Weibel \cite{CHWW2} gave explicit
formulas for $NK_i(R)=K_i(R[X])/K_i(R)$, $i\ge0$, in terms of
Hochschild homology and the cohomology of K\"ahler differentials
in the $cdh$-topology. As a result, they gave affirmative answer
to Bass' question whether $NK_i(R)=0$ implies $R$ is $K$-regular
\cite{Ba2}, when $R$ is an algebra over a large field of
characteristic 0:

\begin{theorem}(\cite[Theorem 0.1(b)]{CHWW2})\label{BassQ}
Let $\kk$ be a characteristic 0 field of infinite transcendence
degree and $R$ an algebra of essentially finite type over $\kk$.
Then we have
$$
NK_i(R)=0\quad\Longrightarrow\quad
K_i(R)=K_i(R[X_1,\ldots,X_n])=0,\ n>0.
$$
\end{theorem}
(An algebra $R$ over a field $\kk$ is of essentially finite type
if it is a localization of a finitely generated algebra over
$\kk$.)

In view of \cite{CHWW1}, for rings as in Theorem \ref{BassQ} and
every index $i\ge\dim R+1$ we obtain the equivalence:
\begin{align*}
NK_i(R)=0\quad \Longleftrightarrow\quad R\ \text{is regular.}
\end{align*}

It is also shown in \cite[Theorem 0.1(a)]{CHWW2} that Theorem
\ref{BassQ} fails for general normal surface singularities over
all number fields, even for $K_0$.

Here we give an alternative short proof of Theorem \ref{BassQ} in
the special case when $R$ is a finitely generated algebra over
$\kk$. The technique we use predates the 1990s.

All our rings are commutative. We also fix the $K$-theoretical
index $i$. The bar over a field refers to the algebraic closure of
the field.

For a functor $\textsf F:\textsf{Rings}\to\textsf{Ab.Groups}$ the
functor $N\F$ is defined by
$$
N\F(R)=\Coker\big(\F(R)\to\F(R[X])\big).
$$
In the notation of \cite[Chap. XII, \S7]{Ba1}, for any ring $R$ we
have natural isomorphisms:
\begin{equation}\label{Iteration}
\F(R[X_1,\ldots,X_n])=(1+N)^n\F(R),\qquad n>0.
\end{equation}

For a ring $R$ the Bloch-Stienstra action \cite{Bl,St} of the ring
of big Witt vectors $\W(R)$ on $NK_i(R)$ makes the latter a
(continuous) $\W(R)$-module. For a systematic theory the reader is
referred to \cite{St}. Weibel \cite{W} extended these actions to a
(continuous) $\W(R)$-module structure on $K_i(A)/K_i(A_0)$ for a
graded $R$-algebra $A=A_0\oplus A_1\oplus\cdots$ with $R\subset
A_0$.

We need the following consequences of the big Witt vectors'
actions:
\begin{proposition}\label{Wactions} Let $n$ be a natural number
and the rings below contain $\QQ$.
\begin{enumerate}
\item[(a)] The groups $K_i(R[X_1\ldots,X_n])/K_i(R)$ are modules
over $R$ and, hence, $\QQ$. \item[(b)] For a Galois extension of
fields $\kk_1\subset\kk_2$ and a noetherian $\kk_1$-algebra $R$
there is a natural $\kk_2$-linear isomorphism:
\begin{align*}
\kk_2\otimes_{\kk_1}\big(K_i(R[X_1,\ldots,&X_n])/K_i(R)\big)=\\
&K_i\big((\kk_2\otimes_{\kk_1}R)[X_1,\ldots,X_n]\big)/K_i\big(\kk_2\otimes_{\kk_1}R\big).
\end{align*}
\item[(c)] To prove Theorem \ref{BassQ} for a noetherian algebra
$R$ over a field $\kk$, without loss of generality we can assume
that $\kk$ is algebraically closed.
\end{enumerate}
\end{proposition}

\begin{proof}
(a) follows from the diagonal embedding $R\to\prod_1^\infty R$,
followed by the ghost isomorphism $\prod_1^\infty R\cong\W(R)$

For (b) we need Galois descent in rational $K$-theory: for a
finite Galois extension of noetherian rings $A\subset B$ (not
necessarily $\QQ\subset A$) with Galois group $G$ we have
\begin{equation}\label{Ratiodescent}
K_p(A)\otimes\QQ=H^0(G,K_p(B)\otimes\QQ).
\end{equation}
This is a special case of Thomason's \'etale descent \cite[Lemma
2.13]{T},\cite[Proposition 11.10]{TT}, worked out in detail for
the reader's convenience in \cite[Lemma 8.4]{G1}.

Returning to our situation, the map
$$
\phi:K_i\big(R[X_1,\ldots,X_n]\big)/K_i(R)\to
K_i\big((\kk_2\otimes_{\kk_1}R)[X_1,\ldots,X_n]\big)/K_i\big(\kk_2\otimes_{\kk_1}R\big)
$$
is $\W(\kk_1)$-linear (compare with formula (1) in \cite{G2}).
Also, if $G$ is the Galois group of $\kk_2$ over $\kk_1$ then
$H^0\big(G,\W(\kk_2)\big)=\W(\kk_1)$. In particular, the target
group of $\phi$ is a Galois $\kk_1$-module. Therefore,
(\ref{Ratiodescent}) implies (b).

(c) follows from (b) because $K_i$ commutes with inductive limits
and so
$$
\bar\kk\otimes_{\kk}\big(K_i(R[X_1,\ldots,X_n])/K_i(R)\big)=
K_i\big((\bar\kk\otimes_{\kk}R)[X_1,\ldots,X_n]\big)/K_i\big(\bar\kk\otimes_{\kk}R\big).
$$
\end{proof}

\begin{remark}\label{Elocal}
Alternatively, Proposition \ref{Wactions}(b) can be derived from
van der Kallen's result \cite[Theorem 3.2]{K} that for an \'etale
extension of rings $A\subset B$, containing $\QQ$, there is an
isomorphism of $\W(B)$-modules
$$
NK_i(B)=\W(B)\otimes_{\W(A)}NK_i(A).
$$
See the comment after formula (5) in \cite{G2}. If one applies van
der Kallen's formula to the situation of Proposition
\ref{Wactions}(b), the difference between the two operator rings
$\kk_2\otimes\W(R)$ and $\W(\kk_2\otimes R)$, resulting from the
two different approaches, is absorbed by the continuity property
of the $NK_i$-groups as modules over the big Witt vectors.
\end{remark}

What follows from this point on is a direct $K_i$-analogue of the
corresponding facts from \cite[\S\S1,2]{MP}:

\begin{lemma}[Affine Horrocks for $K_i$]\label{Monic}
Let $R$ be a ring and $f\in R[X]$ be a monic polynomial with $\deg
f>0$. Then the map $K_i(R[X])\to K_i(R[X]_f)$ is injective.
\end{lemma}

\begin{proof}
\emph{Step 1.} First we prove that the map $K_i(R)\to K_i(R[X]_f)$
is injective.

The rings $R[X]/(f-1)$ and $R[X]/(Xf-1)$ are free $R$-modules of
ranks $\deg f$ and $1+\deg f$, respectively. Therefore, the
transfer maps for $K_i$, corresponding to the two composite
homomorphisms:
$$
\xymatrix{ R\ar[r]& R[X]_f\ar[rr]^{1/f
\mapsto1\qquad}&&R[X]/(f-1),\\ R\ar[r]& R[X]_f\ar[rr]^{1/f
\mapsto\bar X\qquad}&&R[X]/(Xf-1), }
$$
show that $\ker\big(K_i(R)\to K_i(R[X_f])\big)$ is a torsion group
for both exponents $\deg f$ and $1+\deg f$. Therefore, it is a
trivial group.

\medskip\noindent\emph{Step 2.} Here we observe that for a
ring $R$ and two comaximal elements $a,b\in R$ the map
$$
\ker\big(K_i(R)\to
K_i(R_{ab})\big)\longrightarrow\ker\big(K_i(R_a)\to
K_i(R_{ab})\big)
$$
is surjective. This follows directly from Thomason's
Mayer-Vietoris sequence \cite[Theorem 8.1]{TT}, applied to the
open affine cover $\Spec(R_a)\cup\Spec(R_b)=\Spec(R)$:
\begin{align*}
\cdots\to K_i(R)\to K_i(R_a)\oplus K_i(R_b)\to K_i(R_{ab})\to
K_{i-1}(R)\to\cdots
\end{align*}

\medskip\noindent\emph{Notice.} In the proof
of Theorem \ref{BassQ} for finitely generated $\kk$-algebras we
will have $a$ and $b$ both non-zerodivisors.  It was observed by
Murthy \cite{JS} that in such a situation the Mayer-Vietoris
sequence above follows from Quillen's localization sequence.

\medskip\noindent\emph{Step 3.} Assume $\alpha\in\ker\big(K_i(R[X])\to
K_i(R[X]_f)\big)$. Let $\alpha'$ be the image of $\alpha$ under
the map $K_i(R[X])\to K_i(R[X,X^{-1}])$. Then
\begin{align*}
\alpha'\in&\ker\big(K_i(R[X,X^{-1}])\to K_i(R[X]_{Xf})\big)=\\
&\ker\big(K_i(R[X,X^{-1}])\to K_i(R[X^{-1}]_{X^{-1}g})\big),
\end{align*}
where $g=X^{-\deg f}f\in R[X^{-1}]$. Since $f$ is monic (in $X$),
the elements $X^{-1},g\in R[X^{-1}]$ are comaximal. Therefore, by
step 2 there exists $\beta\in K_i(R[X^{-1}])$ such that
$\beta'=\alpha'$, where $\beta'$ is the image of $\beta$ under the
map $K_i(R[X^{-1}])\to K_i(R[X,X^{-1}])$. Then $\alpha\in
K_i(R[X])\cap K_i(R[X^{-1}])=K_i(R)$, the intersection being
considered in $K_i(R[X,X^{-1}])=K_i(R)\oplus NK_i(R)\oplus
NK^-_i(R)\oplus K_{i-1}(R)$ (Fundamental Theorem of $K$-theory
\cite{Gr}). But then Step 1 implies $\alpha=0$.
\end{proof}

\begin{corollary}\label{Injective}
Let $R$ be a $\kk$-algebra, $\chara\kk=0$. Then the maps
$$
N^nK_i(R[X])\to N^nK_i\big(\overline{\kk(X)}\otimes_\kk
R\big),\qquad n>0,
$$
are injective
\end{corollary}

\begin{proof}
Because of the transfer maps for $K_i$ and Proposition
\ref{Wactions}(a), the maps
\begin{align*}
N^nK_i\big(\kk(X)\otimes_\kk R\big)\to
\lim_\to\ &N^nK_i\big(\mathbb K\otimes_\kk R\big)=\\
&N^nK_i\big(\overline{\kk(X)}\otimes_\kk R\big)
\end{align*}
are injective for all $n>0$, where the field $\mathbb K$ in the
inductive limit runs over the finite extensions of $\kk(X)$. So it
is enough to show that the maps
$$
N^nK_i(R[X])\to N^nK_i\big(\kk(X)\otimes_\kk R\big),\qquad n>0,
$$
are injective. But this follows from Lemma \ref{Monic}, applied to
the ambient groups
$$
N^nK_i(R[X])\subset K_i(R[X,X_1,\ldots,X_n]),
$$
which is possible because
$$
\kk(X)\otimes_\kk R[X_1,\ldots,X_n]=\lim_\to
R[X,X_1,\ldots,X_n]_f,
$$
the inductive limit being taken over the ring homomorphisms
$$
R[X,X_1,\ldots,X_n]_f\to R[X,X_1,\ldots,X_n]_{fg},\quad
f,g\in\kk[X].
$$
\end{proof}

\begin{proof}[Proof of Theorem \ref{BassQ} for finitely generated $\kk$-algebras]
By Proposition \ref{Wactions}(c) we can assume
$\kk=\overline{\QQ\big(X_\lambda\big)}_\Lambda$ for an infinite
set of variables $\{X_\lambda\}_\Lambda$, not containing the
variable $X$. Let $R$ be a finitely generated $\kk$-algebra.

To prove Theorem \ref{BassQ}, by (\ref{Iteration}) it is enough to
show $N^nK_i(R)=0$ for all $n>0$. This we prove by induction on
$n$, the case $n=1$ being the hypothesis in Theorem \ref{BassQ}.
Assume $N^{n-1}K_i(R)=0$. For simplicity, denote $\F=N^{n-1}K_i$.

Fix a representation
$$
R=\frac{\kk[t_1,\ldots,t_m]}{(f_1,\ldots,f_r)}.
$$
There is a finite subset
$\{X_1,\ldots,X_n\}\subset\{X_\lambda\}_\Lambda$ such that the
subfield
$$
\kk_0=\overline{\QQ(X_1,\ldots,X_n)}\subset\kk
$$ contains the
coefficients of $f_1,\ldots,f_r$. Fix a bijection of the form
$$
\phi:\{X_\lambda\}_\Lambda\to\{X_\lambda\}_\Lambda\cup\{X\},\quad\phi(X_i)=X_i,\quad
i=1,\ldots,n.
$$
It gives rise to a $\kk_0$-algebra isomorphism
$$
R\cong\overline{\kk(X)}\otimes_\kk
R=\frac{\overline{\kk\big(X,\{X_\lambda\}_\Lambda\big)}[t_1,\ldots,t_m]}{(f_1,\ldots,f_r)}.
$$
So by Corollary \ref{Injective} we have the sequence of injective
maps:
$$
\F(R)\oplus N\F(R)=\F(R[X])\to\F(\overline{\kk(X)}\otimes_\kk
R)\cong\F(R)=0.
$$
\end{proof}

\end{document}